\documentclass[11pt]{amsart}

\usepackage[a4paper,margin=1.10in]{geometry}
\usepackage{amsmath,amssymb,mathtools}
\usepackage{microtype}
\usepackage[hidelinks]{hyperref}
\usepackage[nameinlink,noabbrev]{cleveref}
\hypersetup{pdftitle={Rigidity of Exact Contractivity in the Complex Polynomial Bohnenblust--Hille Inequality},pdfauthor={Daniel M. Pellegrino and Eduardo V. Teixeira}}
\allowdisplaybreaks[1]
\numberwithin{equation}{section}

\newtheorem{theorem}{Theorem}[section]
\newtheorem{lemma}[theorem]{Lemma}
\theoremstyle{remark}
\newtheorem{remark}[theorem]{Remark}
\theoremstyle{plain}
\newtheorem{mainresult}{Theorem}

\newcommand{\M}{\mathcal M}
\newcommand{\dist}{\operatorname{dist}}

\title[Rigidity of exact polynomial BH contractivity]
{Rigidity of Exact Contractivity in the Complex Polynomial Bohnenblust--Hille Inequality}

\author{Daniel M. Pellegrino}
\address{Departamento de Matem\'atica, Universidade Federal da Para\'iba, Jo\~ao Pessoa, PB, Brazil}
\email{daniel.pellegrino@academico.ufpb.br}

\author{Eduardo V. Teixeira}
\address{Department of Mathematics, Oklahoma State University, Stillwater, OK 74078, USA}
\email{eduardo.teixeira@okstate.edu}

\subjclass[2020]{Primary 46G25; Secondary 42C10, 05D40}
\keywords{Bohnenblust--Hille inequality, exact contractivity, homogeneous polynomials, prescribed monomial supports, rigidity}

\begin{document}

\begin{abstract}
The complex polynomial Bohnenblust--Hille inequality bounds the
coefficient $\ell_{2m/(m+1)}$ norm of a complex $m$-homogeneous polynomial
by its supremum norm, with a constant independent of the number of variables.
For each fixed integer $N\ge3$, we characterize the choices of $N$
degree-$m$ monomials for which the inequality holds with optimal
constant $1$ for every polynomial supported on them, provided $m$
exceeds an explicit threshold depending only on $N$.
This holds precisely when their exponent vectors contain no nontrivial
three-term arithmetic progression; equality for a nonzero polynomial
then holds only for monomials.
We also prove that, in the same degree range, every polynomial with at
most $N$ nonzero coefficients that violates the constant-one inequality
has a unique coefficient of maximal modulus, whose exponent is the center
of such a progression in its support. Exact contractivity also follows
from a comparability condition on the two largest coefficients, with a
degree threshold that necessarily depends on the number of nonzero monomials.
\end{abstract}

\maketitle
\setcounter{tocdepth}{2}
\tableofcontents

\section{Introduction}\label{sec:introduction}

The Bohnenblust--Hille inequality originated in Bohr's problem on the
absolute convergence of Dirichlet series. For a Dirichlet series
$D(s)=\sum_{k=1}^{\infty}a_k k^{-s}$ with finite abscissas of absolute
and uniform convergence, denote these abscissas by $\sigma_a(D)$ and
$\sigma_u(D)$, respectively.
Bohr proved in 1913 that $\sigma_a(D)-\sigma_u(D)\le 1/2$ and asked
whether this bound was sharp \cite{Bohr1913}. In 1931, Bohnenblust and
Hille settled the problem by proving that
$\sup_D(\sigma_a(D)-\sigma_u(D))=1/2$ \cite{BohnenblustHille}.
A central ingredient was their extension of Littlewood's $4/3$-inequality
to multilinear forms \cite{Littlewood1930}. Consider a complex
$m$-homogeneous polynomial $P:\mathbb C^n\to\mathbb C$,
\[
 P(z)=\sum_{|\alpha|=m}a_\alpha z^\alpha,
\]
where $\alpha=(\alpha_1,\ldots,\alpha_n)$ has nonnegative integer
entries, $|\alpha|=\alpha_1+\cdots+\alpha_n$, and
$z^\alpha=z_1^{\alpha_1}\cdots z_n^{\alpha_n}$. The polynomial
Bohnenblust--Hille inequality asserts that
\begin{equation}\label{eq:classical-BH}
 \left(\sum_{|\alpha|=m}|a_\alpha|^{2m/(m+1)}\right)^{(m+1)/(2m)}
 \le D_m\|P\|_\infty,
\end{equation}
where $\mathbb D=\{z\in\mathbb C:|z|<1\}$,
$\|P\|_\infty=\sup_{z\in\mathbb D^n}|P(z)|$, and $D_m$ denotes
the optimal constant independent of $n$. H\"older's inequality gives
\[
 \sum_{|\alpha|=m}|a_\alpha|
 \le D_m\binom{n+m-1}{m}^{\frac{m-1}{2m}}\|P\|_\infty.
\]
Thus bounds for $D_m$ control the Sidon constants of the degree-$m$
monomials, that is, the optimal constants comparing the sum of the
coefficient moduli with the supremum norm. They also enter the Bohr
radius problem: finding the largest $r\in(0,1)$ for which
$\sum_\alpha |c_\alpha|r^{|\alpha|}\le\|f\|_\infty$ holds for every
bounded holomorphic function $f(z)=\sum_\alpha c_\alpha z^\alpha$ on
$\mathbb D^n$. These estimates give the asymptotic order
$\sqrt{(\log n)/n}$ for the Bohr radius as $n\to\infty$
\cite{BayartPellegrinoSeoane,DefantFrerickOrtegaOunaiesSeip}.

The bounds for $D_m$ have been progressively improved. In 2011, Defant, Frerick, Ortega-Cerd\`a, Ouna\"ies and Seip established hypercontractivity, showing that $D_m\le C^m$ for some absolute constant $C>1$ \cite{DefantFrerickOrtegaOunaiesSeip}. In 2014, subexponential growth of the complex polynomial constants was established \cite{BayartPellegrinoSeoane}. In earlier joint work \cite[Remark~5.1 and Appendix~A]{PellegrinoTeixeiraPolynomialGrowth}, we proved that, for every $B>1/2$, there exists $K_B>0$ such that
\begin{equation}\label{eq:classical-upper-history}
 D_m\le K_Bm^B,
 \qquad m\ge1,
\end{equation}
and, more precisely, for some absolute constants $C>0$ and $m_0\ge1$,
\[
 D_m\le \sqrt m\exp\!\bigl(C\sqrt{\log m}\,\log\log m\bigr)
 \qquad (m\ge m_0).
\]
For fixed $m,n\ge1$, let $D_{m,n}$ denote the optimal constant
in \eqref{eq:classical-BH} for polynomials on $\mathbb C^n$.
The same work proves that, for any sequence of positive integers $(n_m)_{m\ge1}$,
\begin{equation}\label{eq:classical-lower-history}
 D_{m,n_m}\longrightarrow2
 \qquad\text{whenever}\qquad
 \frac{n_m}{m}\longrightarrow1.
\end{equation}
Since $D_m=\sup_{n\ge1}D_{m,n}$, we have
$\liminf_{m\to\infty}D_m\ge2$.
In particular, $D_m>1$ for all sufficiently large $m$. The exact
asymptotic behavior of $D_m$ remains open.

A parallel line of research examines the Bohnenblust--Hille inequality
under restrictions on the number of variables occurring in each monomial.
Carando, Defant and Sevilla-Peris obtained polynomial growth in the
degree for a fixed number of active variables
\cite{CarandoDefantSevillaPeris}; subsequent work established bounds
uniform in the degree \cite{MaiaNogueiraPellegrino}.
Asymptotic contractivity in this regime follows, by interpolation
with Parseval's identity, from the
support-sensitive inequalities of Defant, Galicer, Mansilla,
Masty\l o and Muro \cite[Theorem~2.12]{DefantGalicerMansillaMastyloMuro}.
Caro-Montoya, N\'u\~nez-Alarc\'on and Serrano-Rodr\'iguez obtained
an improved quantitative estimate for homogeneous polynomials
\cite{CaroMontoyaNunezSerrano}. Here we instead fix the total number
of monomials.

For a monomial, the coefficient norm equals the supremum norm, so the
Bohnenblust--Hille constant cannot be smaller than $1$. Related
constant-one inequalities are known in the completely bounded setting
\cite{ArunachalamDuttEscuderoPalazuelos} and, for sufficiently large
degrees, for symmetric complex multilinear forms
\cite[Corollary~1.3]{NunezPellegrinoRaposoTeixeiraSymmetric}.

For a nonzero polynomial $P$, the Bohnenblust--Hille ratio
\[
 R_m(P):=\frac{\left(\sum_{|\alpha|=m}|a_\alpha|^{2m/(m+1)}\right)^{(m+1)/(2m)}}{\|P\|_\infty}
\]
was examined under normalized surface measure on Euclidean
coefficient spheres in \cite{PellegrinoTeixeiraTypical}. For each
prescribed sequence of monomial supports, $R_m(P)\le1$ in all
sufficiently large degrees for almost every sequence of coefficient
directions with respect to the product of these measures.

We ask when the inequality holds for every choice of coefficients
on a prescribed support. Fix an integer
$N\ge3$. For all sufficiently large degrees $m$, we prove that the
optimal Bohnenblust--Hille constant for polynomials supported on $N$
prescribed degree-$m$ monomials is exactly $1$ if and only if their
exponent vectors contain no nontrivial three-term arithmetic progression.
The degree threshold depends only on $N$, not on the number of variables.

Write
\[
 \mathbb N:=\{1,2,\ldots\},\qquad
 \mathbb N_0:=\mathbb N\cup\{0\},\qquad
 \mathbb T:=\{z\in\mathbb C:|z|=1\}.
\]
For $m,n\in\mathbb N$, put
\begin{equation}\label{eq:Mmn}
 \M_{m,n}:=\{\alpha\in\mathbb N_0^n:|\alpha|=m\},
\end{equation}
and let $\mathcal P_m(\mathbb C^n)$ denote the complex $m$-homogeneous polynomials on $\mathbb C^n$. Every $P\in\mathcal P_m(\mathbb C^n)$ has a unique representation
\begin{equation}\label{eq:polynomial-intro}
 P(z)=\sum_{\alpha\in\M_{m,n}}a_\alpha z^\alpha,
 \qquad z\in\mathbb C^n.
\end{equation}
Define the exponent support of $P$ and, for $1\le p<\infty$,
its coefficient norm by
\begin{equation}\label{eq:coefficient-norm-support}
\begin{aligned}
 \operatorname{supp}P&:=\{\alpha\in\M_{m,n}:a_\alpha\ne0\},\\
 \|P\|_{\ell_p}&:=
 \left(\sum_{\alpha\in\operatorname{supp}P}|a_\alpha|^p\right)^{1/p}.
\end{aligned}
\end{equation}
Recall that the supremum norm satisfies
\begin{equation}\label{eq:supnorm-definition}
 \|P\|_\infty=\sup_{z\in\mathbb D^n}|P(z)|
 =\max_{z\in\mathbb T^n}|P(z)|.
\end{equation}
Orthogonality of monomials on $\mathbb T^n$ gives
\begin{equation}\label{eq:parseval-intro}
 \|P\|_{\ell_2}\le\|P\|_\infty.
\end{equation}

For any nonempty set $\Lambda\subset\M_{m,n}$ of exponent vectors, define
\begin{equation}\label{eq:restricted-D}
 D_m(\Lambda)
 :=\sup\left\{
 \frac{\|P\|_{\ell_{2m/(m+1)}}}{\|P\|_\infty}:
 0\ne P\in\mathcal P_m(\mathbb C^n),\ \operatorname{supp}P\subseteq\Lambda
 \right\}.
\end{equation}

The monomial $P(z)=z^\alpha$, $\alpha\in\Lambda$, shows that
$D_m(\Lambda)\ge1$. We say that the inequality is exactly contractive
on $\Lambda$ when $D_m(\Lambda)=1$.

If $|\Lambda|=N$ and $\operatorname{supp}P\subseteq\Lambda$,
the comparison of coefficient norms and Parseval's identity give
\[
 \|P\|_{\ell_{\frac{2m}{m+1}}}\le N^{1/(2m)}\|P\|_{\ell_2}
 \le N^{1/(2m)}\|P\|_\infty,
 \qquad
 1\le D_m(\Lambda)\le N^{1/(2m)}.
\]
For fixed $N$, the convergence to $1$ is uniform in $n$ and $\Lambda$.
We determine when $D_m(\Lambda)=1$ from the arithmetic of the exponent set.

We use the usual coordinatewise arithmetic on $\mathbb Z^n$.
A nontrivial three-term arithmetic progression in $\Lambda$ is a triple
\begin{equation}\label{eq:3AP}
 \alpha-\beta,\qquad \alpha,\qquad \alpha+\beta\in\Lambda,
 \qquad \beta\in\mathbb Z^n\setminus\{0\}.
\end{equation}
We call $\alpha$ its \emph{center} and $\Lambda$ \emph{progression-free} if it contains no such triple.

\begin{mainresult}\label{thm:A}
For every integer $N\ge3$ there exists an explicit integer
$M_N\ge3$, depending only on $N$, such that for every $m\ge M_N$,
$n\ge2$, and every set $\Lambda\subset\M_{m,n}$ with
$|\Lambda|=N$, the following assertions hold.
\begin{enumerate}
\item[\textup{(1)}] If $\Lambda$ contains no nontrivial three-term arithmetic progression, every nonzero $m$-homogeneous polynomial $P:\mathbb C^n\to\mathbb C$ with $\operatorname{supp}P\subseteq\Lambda$ satisfies
\begin{equation}\label{eq:theorem-contractivity}
 \left(\sum_{|\alpha|=m}|a_\alpha|^{\frac{2m}{m+1}}\right)^{\frac{m+1}{2m}}
 \le\|P\|_\infty.
\end{equation}
Equality holds if and only if $P$ is a monomial. In particular, $D_m(\Lambda)=1$.

\item[\textup{(2)}] If $\Lambda$ contains a nontrivial three-term arithmetic progression, then $D_m(\Lambda)>1$; that is, some $m$-homogeneous polynomial $P:\mathbb C^n\to\mathbb C$ with $\operatorname{supp}P\subseteq\Lambda$ violates \eqref{eq:theorem-contractivity}.

\item[\textup{(3)}] Every polynomial $P(z)=\sum_{\gamma\in\Lambda}a_\gamma z^\gamma$ violating \eqref{eq:theorem-contractivity} has a unique coefficient of largest modulus, and its associated exponent is the center of a nontrivial three-term arithmetic progression in $\operatorname{supp}P$.
\end{enumerate}
\end{mainresult}

An admissible threshold $M_N$ is given in \eqref{eq:theta-M}.
The obstruction in part~\textup{(2)} applies in every degree for which such
a progression exists.

\begin{remark}\label{rem:real-obstruction}
The complex field is essential in Theorem~\ref{thm:A}. Indeed, for every even integer $m\ge4$, the set $\Lambda_m=\{(m,0),(m-1,1),(0,m)\}$ contains no nontrivial three-term arithmetic progression, but the real polynomial $P_m(x,y)=x^m-y^m$, with $\operatorname{supp}P_m\subseteq\Lambda_m$, satisfies
\[
 \|P_m\|_{\infty,\mathbb R}
 :=\sup_{(x,y)\in[-1,1]^2}|x^m-y^m|=1
 <2^{\frac{m+1}{2m}}=\|P_m\|_{\ell_{\frac{2m}{m+1}}}.
\]
\end{remark}

A condition on the two largest coefficient moduli gives exact contractivity without assuming progression-freeness.

\begin{mainresult}\label{thm:B}
Let $(M_N)_{N\ge3}$ be the degree thresholds in Theorem~\ref{thm:A}.
\begin{enumerate}
\item[\textup{(1)}] Fix an integer $N\ge3$. Let $m\ge M_N$,
$n\ge1$, and let
$P:\mathbb C^n\to\mathbb C$ be a nonzero $m$-homogeneous
polynomial with at most $N$ nonzero coefficients. If $P$ is a monomial,
or if its two largest coefficient moduli, $|a_1|\ge|a_2|>0$, satisfy
\begin{equation}\label{eq:theorem-B-condition}
 |a_1|\le 4(N-1)^{3/2}|a_2|,
\end{equation}
then
\begin{equation}\label{eq:theorem-B-contractivity}
 \left(\sum_{|\alpha|=m}|a_\alpha|^{\frac{2m}{m+1}}\right)^{\frac{m+1}{2m}}
 \le\|P\|_\infty.
\end{equation}
In particular, the inequality holds whenever the maximum coefficient
modulus is attained at least twice. The constant $1$ is optimal for
this class.

\item[\textup{(2)}] The degree threshold in \textup{(1)} cannot be chosen
independently of $N$: for every integer $M\ge1$ there exist $m\ge M$,
$n\ge2$, and a complex $m$-homogeneous polynomial
$P:\mathbb C^n\to\mathbb C$ whose largest coefficient modulus is
attained at least twice, but
\begin{equation}\label{eq:theorem-B-nonuniform}
 \left(\sum_{|\alpha|=m}|a_\alpha|^{\frac{2m}{m+1}}\right)^{\frac{m+1}{2m}}
 >\|P\|_\infty.
\end{equation}
\end{enumerate}
\end{mainresult}

\section{Auxiliary results}\label{sec:auxiliary}

For a fixed number of coefficients, near equality in Parseval's
inequality forces concentration on a single monomial. After factoring
out that monomial, a three-term progression centered at its exponent
becomes a pair of opposite
frequencies. Their coefficients can be chosen so that the perturbation is
purely imaginary on the torus: its contribution to the squared supremum
norm is then quadratic, while its contribution to the coefficient norm has
order $2m/(m+1)<2$. The reflection estimates below show that, in the
concentrated regime, a violation requires such a pair.

\begin{lemma}\label{lem:progression-witness}
Let $\Lambda\subset\M_{m,n}$ contain
$\alpha-\beta,\alpha,\alpha+\beta$, where $\beta\in\mathbb Z^n\setminus\{0\}$.
Then there is an $m$-homogeneous polynomial $P$ supported on these three
monomials such that
\begin{equation}\label{eq:exact-witness}
 \frac{\|P\|_{\ell_{\frac{2m}{m+1}}}}{\|P\|_\infty}
 =\big(1+2^{1-m}\big)^{1/(2m)}>1.
\end{equation}
Moreover, for every fixed $1\le p<2$, there is a polynomial supported on the same three monomials for which
\[
 \|P\|_{\ell_p}>\|P\|_\infty.
\]
\end{lemma}

\begin{proof}
For $\varepsilon>0$, set
\[
 P_\varepsilon(z)=z^\alpha+\varepsilon z^{\alpha+\beta}
 -\varepsilon z^{\alpha-\beta}.
\]
Since $\beta\ne0$, the values of $z^\beta$ on $\mathbb T^n$
fill $\mathbb T$. Hence
\[
\begin{aligned}
 \|P_\varepsilon\|_\infty^2
 &=\max_{z\in\mathbb T^n}
   |z^\alpha|^2\big|1+\varepsilon z^\beta
                           -\varepsilon z^{-\beta}\big|^2\\
 &=\max_{|w|=1}\big|1+\varepsilon(w-w^{-1})\big|^2\\
 &=\max_{\vartheta\in\mathbb R}
   \big|1+2i\varepsilon\sin\vartheta\big|^2\\
 &=\max_{\vartheta\in\mathbb R}
   \big(1+4\varepsilon^2\sin^2\vartheta\big)
  =1+4\varepsilon^2.
\end{aligned}
\]
Choosing $\varepsilon^2=2^{-(m+1)}$ and writing
$q_m:=\frac{2m}{m+1}$, we obtain
\[
 2\varepsilon^{q_m}=4\varepsilon^2=2^{1-m},
\]
and hence
\[
\frac{\|P_\varepsilon\|_{\ell_{q_m}}^2}{\|P_\varepsilon\|_\infty^2}
=\frac{(1+2\varepsilon^{q_m})^{2/q_m}}{1+4\varepsilon^2}
=\big(1+2^{1-m}\big)^{1/m}>1.
\]

Fix $1\le p<2$. Bernoulli's inequality,
$(1+x)^{2/p}\ge1+(2/p)x$ for $x\ge0$, gives
\begin{equation}\label{eq:p-obstruction}
\begin{aligned}
\frac{\|P_\varepsilon\|_{\ell_p}^2}{\|P_\varepsilon\|_\infty^2}
&=\frac{(1+2\varepsilon^p)^{2/p}}{1+4\varepsilon^2}
\ge\frac{1+(4/p)\varepsilon^p}{1+4\varepsilon^2}\\
&=1+\frac{4\varepsilon^p}{1+4\varepsilon^2}
\left(\frac1p-\varepsilon^{2-p}\right)>1
\end{aligned}
\end{equation}
for $0<\varepsilon<p^{-1/(2-p)}$.
\end{proof}

At an extreme difference of the exponent set, the corresponding Fourier
coefficient of $|P|^2$ contains a single product. Positivity of
$\|P\|_\infty^2-|P|^2$ then gives the following estimate, adapted from
\cite[Lemma~3.8]{PellegrinoTeixeiraTypical}.

\begin{lemma}\label{lem:exposed-pair}
Let $\Lambda\subset\M_{m,n}$ have at least two elements. There are distinct
$\alpha_+,\alpha_-\in\Lambda$, depending only on $\Lambda$, such that every
\[
 P(z)=\sum_{\alpha\in\Lambda}a_\alpha z^\alpha
\]
satisfies
\begin{equation}\label{eq:exposed-pair}
 \|P\|_\infty^2
 \ge \|P\|_{\ell_2}^2+|a_{\alpha_+}a_{\alpha_-}|.
\end{equation}
\end{lemma}

\begin{proof}
Order the elements of $\Lambda$ lexicographically: compare first coordinates,
then second coordinates if the first agree, and continue in this way. Let
$\alpha_-$ be the first exponent and $\alpha_+$ the last, and put
$\delta=\alpha_+-\alpha_-\ne0$.

If $\alpha,\gamma\in\Lambda$ and $\alpha-\gamma=\delta$, then
\[
 \alpha_+-\alpha=\alpha_--\gamma.
\]
In lexicographic order the left-hand side is nonnegative and the right-hand
side is nonpositive. They must therefore both vanish. Thus the difference
$\delta$ is represented only by the pair
$(\alpha,\gamma)=(\alpha_+,\alpha_-)$.

Let $\mu_n$ denote normalized Haar probability measure on
$\mathbb T^n$. For $\xi\in\mathbb Z^n$,
\[
 \int_{\mathbb T^n}z^\xi\,d\mu_n(z)
 =\prod_{j=1}^n\int_{\mathbb T}z_j^{\xi_j}\,d\mu_1(z_j)
 =\begin{cases}1,&\xi=0,\\0,&\xi\ne0.\end{cases}
\]
By Fourier orthogonality and the uniqueness of the pair
$(\alpha_+,\alpha_-)$,
\begin{equation}\label{eq:exposed-fourier}
\begin{aligned}
 \int_{\mathbb T^n}|P(z)|^2z^{-\delta}\,d\mu_n(z)
 &=\sum_{\alpha,\gamma\in\Lambda}a_\alpha\overline{a_\gamma}
   \int_{\mathbb T^n}z^{\alpha-\gamma-\delta}\,d\mu_n(z)\\
 &=\sum_{\substack{\alpha,\gamma\in\Lambda\\\alpha-\gamma=\delta}}
   a_\alpha\overline{a_\gamma}
 =a_{\alpha_+}\overline{a_{\alpha_-}}.
\end{aligned}
\end{equation}
Since $\delta\ne0$, orthogonality gives
\begin{equation}\label{eq:zero-character}
\int_{\mathbb T^n}z^{-\delta}\,d\mu_n(z)=0.
\end{equation}
By \eqref{eq:exposed-fourier} and \eqref{eq:zero-character}, and since
$\|P\|_\infty^2-|P(z)|^2\ge0$ and $|z^{-\delta}|=1$ on $\mathbb T^n$,
\begin{align*}
 |a_{\alpha_+}a_{\alpha_-}|
 &=\left|\int_{\mathbb T^n}
    \bigl(\|P\|_\infty^2-|P(z)|^2\bigr)
    z^{-\delta}\,d\mu_n(z)\right|\\
 &\le\int_{\mathbb T^n}
    \bigl(\|P\|_\infty^2-|P(z)|^2\bigr)\,d\mu_n(z)\\
 &=\|P\|_\infty^2-\int_{\mathbb T^n}|P(z)|^2\,d\mu_n(z)\\
 &=\|P\|_\infty^2-\sum_{\alpha\in\Lambda}|a_\alpha|^2
 =\|P\|_\infty^2-\|P\|_{\ell_2}^2.
\end{align*}
\end{proof}

\begin{lemma}\label{lem:concentration}
Let $\Lambda\subset\M_{m,n}$ have cardinality $N\ge2$, and let
the polynomial
\[
 P(z)=\sum_{\alpha\in\Lambda}a_\alpha z^\alpha
\]
satisfy $\|P\|_{\ell_2}=1$.
\begin{enumerate}
\item[\textup{(1)}] The coefficients of $P$ satisfy
\begin{equation}\label{eq:concentration}
 1-\max_{\alpha\in\Lambda}|a_\alpha|^2
 \le 25N(N-1)\bigl(\|P\|_\infty^2-1\bigr)^{2^{-(N-2)}}.
\end{equation}

\item[\textup{(2)}] If $m\ge2$ and
\[
 \|P\|_{\ell_{\frac{2m}{m+1}}}>\|P\|_\infty,
\]
then
\begin{equation}\label{eq:exceptional-concentration}
 1-\max_{\alpha\in\Lambda}|a_\alpha|^2
 <25N(N-1)\bigl(N^{1/m}-1\bigr)^{2^{-(N-2)}}.
\end{equation}
\end{enumerate}
\end{lemma}

\begin{proof}
\textup{(1)} 
By Parseval's identity,
$\Delta:=\|P\|_\infty^2-1\ge0$. If $\Delta\ge1$, then
\[
 1-\max_{\alpha\in\Lambda}|a_\alpha|^2
 \le1\le25N(N-1)\Delta^{2^{-(N-2)}}.
\]
Hence we may assume $0\le\Delta<1$.

We apply Lemma~\ref{lem:exposed-pair} repeatedly, each time deleting
the coefficient of smaller modulus from the selected pair.
Let $\Lambda_0=\Lambda$ and $Q_0=P$. For $0\le k\le N-2$,
let $\alpha_{k,-}$ and $\alpha_{k,+}$ be the first and last elements
of $\Lambda_k$ in lexicographic order, and put
\[
 \eta_k:=
 \begin{cases}
  \alpha_{k,+},& |a_{\alpha_{k,+}}|\le |a_{\alpha_{k,-}}|,\\
  \alpha_{k,-},& |a_{\alpha_{k,+}}|>|a_{\alpha_{k,-}}|.
 \end{cases}
\]
Define
\[
 \Lambda_{k+1}=\Lambda_k\setminus\{\eta_k\},\qquad
 c_k=a_{\eta_k},\qquad
 Q_{k+1}=Q_k-c_k z^{\eta_k}.
\]
Each step removes one exponent, so $|\Lambda_k|=N-k$ and
\[
 Q_k(z)=\sum_{\alpha\in\Lambda_k}a_\alpha z^\alpha,
 \qquad \|Q_k\|_{\ell_2}^2
 =\sum_{\alpha\in\Lambda_k}|a_\alpha|^2\le1.
\]
Cauchy--Schwarz gives
\begin{equation}\label{eq:Qk-CS}
 \|Q_k\|_\infty\le\|Q_k\|_{\ell_1}
 \le\sqrt{N-k}\,\|Q_k\|_{\ell_2}\le\sqrt N,
 \qquad 0\le k\le N-1.
\end{equation}

For $0\le k\le N-1$, put
\[
 d_k:=\|Q_k\|_\infty^2-\|Q_k\|_{\ell_2}^2.
\]
Parseval's identity and \eqref{eq:Qk-CS} give
$0\le d_k\le N$, and $d_0=\Delta$.
For the pair selected above, Lemma~\ref{lem:exposed-pair} gives
\begin{equation}\label{eq:deletion-coefficient}
\begin{aligned}
 |c_k|^2
 &=\min\bigl\{|a_{\alpha_{k,+}}|^2,|a_{\alpha_{k,-}}|^2\bigr\}\\
 &\le|a_{\alpha_{k,+}}a_{\alpha_{k,-}}|
 \le\|Q_k\|_\infty^2-\|Q_k\|_{\ell_2}^2=d_k,
 \quad 0\le k\le N-2.
\end{aligned}
\end{equation}

Since $Q_{k+1}=Q_k-c_kz^{\eta_k}$, we have
$\|Q_{k+1}\|_\infty\le\|Q_k\|_\infty+|c_k|$ and
$\|Q_{k+1}\|_{\ell_2}^2=\|Q_k\|_{\ell_2}^2-|c_k|^2$.
Using \eqref{eq:Qk-CS} and \eqref{eq:deletion-coefficient}, for $0\le k\le N-2$,
\begin{align*}
 d_{k+1}
 &\le (\|Q_k\|_\infty+|c_k|)^2
       -(\|Q_k\|_{\ell_2}^2-|c_k|^2)\\
 &=d_k+2\|Q_k\|_\infty|c_k|+2|c_k|^2\\
 &\le3d_k+2\sqrt N\sqrt{d_k}\\
 &\le 5\sqrt N\sqrt{d_k},
\end{align*}
where the last inequality uses $d_k\le N$.

We prove by induction that, for $0\le k\le N-1$,
\begin{equation}\label{eq:dk}
 d_k\le(5\sqrt N)^{2(1-2^{-k})}\Delta^{2^{-k}}
 \le25N\Delta^{2^{-k}}.
\end{equation}
For $k=0$, the first inequality is $d_0=\Delta$.
If the first inequality holds for some $0\le k\le N-2$, then
\begin{align*}
 d_{k+1}
 &\le5\sqrt N\sqrt{d_k}\\
 &\le5\sqrt N
 \bigl((5\sqrt N)^{2(1-2^{-k})}\Delta^{2^{-k}}\bigr)^{1/2}\\
 &=(5\sqrt N)^{2(1-2^{-(k+1)})}\Delta^{2^{-(k+1)}}.
\end{align*}
The second inequality in \eqref{eq:dk} follows from
$2(1-2^{-k})\le2$ and $5\sqrt N\ge1$.

The set $\Lambda_{N-1}$ has one element, say $\eta_*$. By the construction of
$Q_k$,
\[
 Q_{N-1}(z)=a_{\eta_*}z^{\eta_*},\qquad
 P(z)=a_{\eta_*}z^{\eta_*}
      +\sum_{k=0}^{N-2}c_kz^{\eta_k}.
\]
The exponents $\eta_0,\ldots,\eta_{N-2},\eta_*$ are distinct. Since
$\|P\|_{\ell_2}=1$, it follows that
\[
 1=|a_{\eta_*}|^2+\sum_{k=0}^{N-2}|c_k|^2.
\]
Since $0\le\Delta<1$, for $0\le k\le N-2$ we have
$\Delta^{2^{-k}}\le\Delta^{2^{-(N-2)}}$. Thus
\begin{align*}
 1-\max_{\alpha\in\Lambda}|a_\alpha|^2
 &\le1-|a_{\eta_*}|^2
 =\sum_{k=0}^{N-2}|c_k|^2\\
 &\le\sum_{k=0}^{N-2}d_k
 \le25N\sum_{k=0}^{N-2}\Delta^{2^{-k}}\\
 &\le25N(N-1)\Delta^{2^{-(N-2)}}.
\end{align*}

\textup{(2)} 
Put $q_m=2m/(m+1)$. H\"older's inequality gives
\[
 \|P\|_{\ell_{q_m}}
 \le N^{1/(2m)}\|P\|_{\ell_2}=N^{1/(2m)}.
\]
Since $\|P\|_{\ell_{q_m}}>\|P\|_\infty$, we obtain
\[
 0\le\|P\|_\infty^2-1
 <\|P\|_{\ell_{q_m}}^2-1\le N^{1/m}-1.
\]
Therefore, by (1),
\begin{align*}
 1-\max_{\alpha\in\Lambda}|a_\alpha|^2
 &\le25N(N-1)\bigl(\|P\|_\infty^2-1\bigr)^{2^{-(N-2)}}\\
 &<25N(N-1)\bigl(N^{1/m}-1\bigr)^{2^{-(N-2)}}.
\end{align*}
\end{proof}

\begin{lemma}\label{lem:reflection}
Let $m\ge2$, let $\Lambda\subset\M_{m,n}$ have $N\ge3$ elements, and let
$P(z)=\sum_{\gamma\in\Lambda}a_\gamma z^\gamma$ satisfy
$\|P\|_{\ell_{\frac{2m}{m+1}}}\ge\|P\|_\infty$.
Fix $\alpha\in\Lambda$ with $a_\alpha\ne0$, and put
\[
 I_\alpha=\{\beta\in\mathbb Z^n\setminus\{0\}:\alpha+\beta\in\Lambda\},
 \qquad b_\beta=\frac{a_{\alpha+\beta}}{a_\alpha}\quad(\beta\in I_\alpha),
\]
with $b_\beta=0$ for $\beta\notin I_\alpha$. Set
\[
 s=\frac{m}{m+1},\qquad
 u=\sum_{\beta\in I_\alpha}|b_\beta|^{2s},\qquad
 \eta=(1+u)^{1/s}-1.
\]
Then, for every $\beta\in I_\alpha$,
\begin{equation}\label{eq:reflection-coefficient}
 |b_\beta+\overline{b_{-\beta}}|\le\eta.
\end{equation}
\end{lemma}

\begin{proof}
Define $F:\mathbb T^n\to\mathbb C$ by
\[
 F(z)=\sum_{\beta\in I_\alpha}b_\beta z^\beta,
 \qquad \frac{P(z)}{a_\alpha z^\alpha}=1+F(z).
\]
Put $q_m=2m/(m+1)=2s$. Since $|z^\alpha|=1$ on
$\mathbb T^n$ and $P(z)=a_\alpha z^\alpha(1+F(z))$,
\[
\begin{aligned}
1+2\operatorname{Re}F(z)+|F(z)|^2
 &=|1+F(z)|^2=\frac{|P(z)|^2}{|a_\alpha|^2}\\
 &\le\frac{\|P\|_\infty^2}{|a_\alpha|^2}
 \le\frac{\|P\|_{\ell_{q_m}}^2}{|a_\alpha|^2}\\
 &=\left(1+\sum_{\beta\in I_\alpha}|b_\beta|^{q_m}\right)^{2/q_m}
 = (1+u)^{1/s}=1+\eta.
\end{aligned}
\]
Consequently,
\[
 H(z):=\frac\eta2-\operatorname{Re}F(z)
 \ge\frac12|F(z)|^2\ge0.
\]
Every $\rho\in I_\alpha$ is nonzero. Orthogonality therefore gives
\[
\int_{\mathbb T^n}F(z)\,d\mu_n(z)
 =\sum_{\rho\in I_\alpha}b_\rho
   \int_{\mathbb T^n}z^\rho\,d\mu_n(z)=0,
 \qquad \int_{\mathbb T^n}H(z)\,d\mu_n(z)=\frac\eta2.
\]
For $\beta\in I_\alpha$, orthogonality gives
\[
\begin{aligned}
\int_{\mathbb T^n}F(z)z^{-\beta}\,d\mu_n(z)
 &=\sum_{\rho\in I_\alpha}b_\rho
   \int_{\mathbb T^n}z^{\rho-\beta}\,d\mu_n(z)=b_\beta,\\
\int_{\mathbb T^n}\overline{F(z)}z^{-\beta}\,d\mu_n(z)
 &=\sum_{\rho\in I_\alpha}\overline{b_\rho}
   \int_{\mathbb T^n}z^{-\rho-\beta}\,d\mu_n(z)=\overline{b_{-\beta}}.
\end{aligned}
\]
Since $\int_{\mathbb T^n}z^{-\beta}\,d\mu_n(z)=0$,
\[
\int_{\mathbb T^n}H(z)z^{-\beta}\,d\mu_n(z)
 =-\frac12\left(b_\beta+\overline{b_{-\beta}}\right).
\]
Since $H\ge0$ and $|z^{-\beta}|=1$ on $\mathbb T^n$,
\[
 |b_\beta+\overline{b_{-\beta}}|
 =2\left|\int_{\mathbb T^n}H(z)z^{-\beta}\,d\mu_n(z)\right|
 \le2\int_{\mathbb T^n}H(z)\,d\mu_n(z)=\eta.
\]
\end{proof}

\begin{lemma}\label{lem:reflection-distance}
Under the hypotheses and notation of Lemma~\ref{lem:reflection}, put
$b=(b_\beta)_{\beta\in I_\alpha}$ and define the real linear subspace
\[
 V_\alpha=\left\{v=(v_\beta)_{\beta\in I_\alpha}\in\mathbb C^{I_\alpha}:
 \begin{array}{ll}
 v_{-\beta}=-\overline{v_\beta},&\beta,-\beta\in I_\alpha,\\
 v_\beta=0,&\beta\in I_\alpha,\ -\beta\notin I_\alpha
 \end{array}
 \right\}.
\]
\begin{enumerate}
\item[\textup{(1)}] We have
\[
 V_\alpha\ne\{0\}
 \quad\Longleftrightarrow\quad
 \alpha\text{ is the center of a nontrivial three-term progression in }\Lambda.
\]
Moreover, with $\dist_2(b,V_\alpha):=\inf_{v\in V_\alpha}\|b-v\|_2$,
\begin{equation}\label{eq:reflection-distance}
 \dist_2(b,V_\alpha)\le\sqrt{N-1}\,\eta.
\end{equation}

\item[\textup{(2)}] Put $t=\sum_{\beta\in I_\alpha}|b_\beta|^2$.
If $m\ge3$, $t\le1$, and $u\le1$, then
\begin{equation}\label{eq:reflection-t}
 \dist_2(b,V_\alpha)\le2(N-1)^{3/4}t^{3/4}.
\end{equation}
\end{enumerate}
\end{lemma}

\begin{proof}
\textup{(1)} 
If $v\in V_\alpha\setminus\{0\}$, some $v_\beta\ne0$,
so $\beta,-\beta\in I_\alpha$. Conversely, if
$\beta,-\beta\in I_\alpha$, the vector with
$v_\beta=1$, $v_{-\beta}=-1$ and all remaining coordinates zero
belongs to $V_\alpha\setminus\{0\}$. The condition
$\beta,-\beta\in I_\alpha$ is equivalent to
$\alpha-\beta,\alpha,\alpha+\beta\in\Lambda$.

For each pair $\{\beta,-\beta\}\subset I_\alpha$ and
$\zeta\in\mathbb C$, we have
\[
 |b_\beta-\zeta|^2+|b_{-\beta}+\overline\zeta|^2
 =2\left|\zeta-\frac{b_\beta-\overline{b_{-\beta}}}{2}\right|^2
 +\frac12|b_\beta+\overline{b_{-\beta}}|^2.
\]
The minimum is attained at
$\zeta=(b_\beta-\overline{b_{-\beta}})/2$ and equals
$\frac12|b_\beta+\overline{b_{-\beta}}|^2$.
Let $\mathcal R$ contain one exponent from each pair
$\{\beta,-\beta\}\subset I_\alpha$, and put
$\mathcal U=\{\beta\in I_\alpha:-\beta\notin I_\alpha\}$.
For $\beta\in\mathcal U$, Lemma~\ref{lem:reflection} gives
$|b_\beta|=|b_\beta+\overline{b_{-\beta}}|\le\eta$.
The minimization separates over the pairs, while every vector in
$V_\alpha$ has zero coordinates on $\mathcal U$. Since
$2|\mathcal R|+|\mathcal U|=|I_\alpha|=N-1$, we obtain
\[
\begin{aligned}
 \dist_2(b,V_\alpha)^2
 &=\frac12\sum_{\beta\in\mathcal R}
     |b_\beta+\overline{b_{-\beta}}|^2
   +\sum_{\beta\in\mathcal U}|b_\beta|^2\\
 &\le\left(\frac{|\mathcal R|}{2}+|\mathcal U|\right)\eta^2
 \le\left(2|\mathcal R|+|\mathcal U|\right)\eta^2\\
 &=(N-1)\eta^2.
\end{aligned}
\]
Taking square roots gives
$\dist_2(b,V_\alpha)\le\sqrt{N-1}\,\eta$.

\textup{(2)} 
Put $J=N-1$. H\"older's inequality, with exponents $1/s$ and
$1/(1-s)$, gives
\[
 u=\sum_{\beta\in I_\alpha}(|b_\beta|^2)^s
 \le\left(\sum_{\beta\in I_\alpha}|b_\beta|^2\right)^s
      \left(\sum_{\beta\in I_\alpha}1\right)^{1-s}
 =t^sJ^{1-s}.
\]
As $3/4\le s<1$, $J\ge1$, and $0\le t\le1$,
$u\le J^{1-s}t^s\le J^{1/4}t^{3/4}$.
For $0\le u\le1$, integration of the derivative of
$x\mapsto(1+x)^{1/s}$ gives
\[
 \eta=(1+u)^{1/s}-1
 =\frac1s\int_0^u(1+x)^{1/s-1}\,dx
 \le\frac43\,2^{1/3}u\le2u,
\]
because $1/s\le4/3$ and $1/s-1\le1/3$. Thus
\[
 \dist_2(b,V_\alpha)
 \le\sqrt J\,\eta\le2\sqrt J\,u\le2J^{3/4}t^{3/4}.
\]
\end{proof}

\begin{lemma}\label{lem:noncenter}
Let $m\ge3$, let $\Lambda\subset\M_{m,n}$ have cardinality $N\ge3$,
and let $P(z)=\sum_{\gamma\in\Lambda}a_\gamma z^\gamma$ satisfy
\[
 \|P\|_{\ell_{\frac{2m}{m+1}}}\ge\|P\|_\infty.
\]
If $\alpha\in\Lambda$ has $a_\alpha\ne0$ and is not the center of a
nontrivial three-term arithmetic progression in $\Lambda$, then either
$P=a_\alpha z^\alpha$ or
\begin{equation}\label{eq:t-lower}
 \frac{\sum_{\gamma\in\Lambda\setminus\{\alpha\}}|a_\gamma|^2}
      {\sum_{\gamma\in\Lambda}|a_\gamma|^2}
 \ge\frac1{16(N-1)^3+1}.
\end{equation}
\end{lemma}

\begin{proof}
Assume $P\ne a_\alpha z^\alpha$. Put $J=N-1$,
$s=m/(m+1)$, and
\[
 t=\frac{\sum_{\gamma\in\Lambda\setminus\{\alpha\}}|a_\gamma|^2}
         {|a_\alpha|^2}>0,
 \qquad
 u=\frac{\sum_{\gamma\in\Lambda\setminus\{\alpha\}}|a_\gamma|^{2s}}
         {|a_\alpha|^{2s}}.
\]
With $I_\alpha$ and $b$ defined in Lemma~\ref{lem:reflection},
and $V_\alpha$ defined in Lemma~\ref{lem:reflection-distance}, the absence of a progression centered at
$\alpha$ gives $V_\alpha=\{0\}$.

If $t<1/(16J^3)$, then $t<1$, and H\"older's inequality gives
\[
\begin{aligned}
 u&\le J^{1-s}t^s\le J^{1/4}t^{3/4}\\
 &<J^{1/4}\left(\frac1{16J^3}\right)^{3/4}
 =\frac1{8J^2}\le1.
\end{aligned}
\]
Since $t>0$, Lemma~\ref{lem:reflection-distance}(2) gives
\[
 1=\frac{\dist_2(b,V_\alpha)}{\sqrt t}
 \le2J^{3/4}t^{1/4}<2J^{3/4}(16J^3)^{-1/4}=1,
\]
a contradiction. Therefore $t\ge1/(16J^3)$, and
\[
 \frac{\sum_{\gamma\in\Lambda\setminus\{\alpha\}}|a_\gamma|^2}
      {\sum_{\gamma\in\Lambda}|a_\gamma|^2}
 =\frac{t}{1+t}\ge\frac1{16J^3+1}
 =\frac1{16(N-1)^3+1}.
\]
\end{proof}

\section{Proof of Theorem A}\label{sec:proofA}

\begin{proof}
Put $q_m:=2m/(m+1)$. For every integer $N\ge3$, put
\begin{equation}\label{eq:constants}
 C_N=25N(N-1),\qquad
 \kappa_N=\frac1{16(N-1)^3+1},
\end{equation}
and define
\begin{equation}\label{eq:theta-M}
 \theta_N=\left(\frac{\kappa_N}{C_N}\right)^{2^{N-2}},\qquad
 M_N=\max\left\{3,\,
 1+\left\lceil\frac{\log N}{\log(1+\theta_N)}\right\rceil\right\}.
\end{equation}
Since $C_N$ increases, $\kappa_N$ decreases, and
$0<\kappa_N/C_N<1$, the sequence $(\theta_N)_{N\ge3}$ decreases.
Thus $\log N/\log(1+\theta_N)$ increases, so $(M_N)_{N\ge3}$
is nondecreasing.
For $m\ge M_N$,
\begin{equation}\label{eq:threshold-use}
 \frac{\log N}{m}<\log(1+\theta_N),
 \qquad N^{1/m}-1<\theta_N<1.
\end{equation}

To prove part~\textup{(3)}, suppose that $P$ violates
\eqref{eq:theorem-contractivity}, and write
$S:=\operatorname{supp}P$, $k:=|S|$. A monomial attains equality.
If $k=2$, write $S=\{\rho,\sigma\}$. Since
$\{z^{\sigma-\rho}:z\in\mathbb T^n\}=\mathbb T$,
\[
 \|P\|_\infty
 =\max_{|w|=1}|a_\rho+a_\sigma w|
 =|a_\rho|+|a_\sigma|
 =\|P\|_{\ell_1}\ge\|P\|_{\ell_{q_m}}.
\]
Thus $3\le k\le N$, and $m\ge M_N\ge M_k$.
Replacing $P$ by $P/\|P\|_{\ell_2}$ preserves the strict inequality,
so we may assume $\sum_{\gamma\in S}|a_\gamma|^2=1$.
Choose $\alpha\in S$ with $|a_\alpha|=\max_{\gamma\in S}|a_\gamma|$,
and set $\tau=1-|a_\alpha|^2$.
The coefficient-norm comparison and the strict violation give
\[
 0\le\Delta:=\|P\|_\infty^2-1
 <\|P\|_{\ell_{q_m}}^2-1
 \le k^{1/m}-1<\theta_k.
\]
Consequently Lemma~\ref{lem:concentration} yields
\begin{equation}\label{eq:tau-upper}
 \tau\le C_k\Delta^{2^{-(k-2)}}
 <C_k\theta_k^{2^{-(k-2)}}=\kappa_k<\frac12.
\end{equation}
For every $\gamma\ne\alpha$ with $\gamma\in S$,
\[
 |a_\gamma|^2\le\tau<\frac12
 <1-\tau=|a_\alpha|^2,
\]
so the largest coefficient modulus is attained only at $\alpha$.
If $\alpha$ were not the center of a three-term progression in $S$,
Lemma~\ref{lem:noncenter}, with $\Lambda=S$, would imply
\[
 \tau=\frac{\sum_{\gamma\in S\setminus\{\alpha\}}|a_\gamma|^2}
          {\sum_{\gamma\in S}|a_\gamma|^2}
 \ge\kappa_k,
\]
contrary to $\tau<\kappa_k$. Thus $\alpha$ is the center of a
nontrivial three-term progression in $\operatorname{supp}P$.

If $\Lambda$ is progression-free, the strict violation is impossible;
hence $D_m(\Lambda)=1$. If $\Lambda$ contains a nontrivial
three-term progression, Lemma~\ref{lem:progression-witness} gives a
polynomial violating \eqref{eq:theorem-contractivity}, and thus
$D_m(\Lambda)>1$.

Assume that $\Lambda$ is progression-free and
$P\ne0$ satisfies $\|P\|_{\ell_{q_m}}=\|P\|_\infty$.
Let $S=\operatorname{supp}P$ and $k=|S|$. If $k=2$, then
\[
 \|P\|_\infty=\|P\|_{\ell_1}>
 \|P\|_{\ell_{q_m}},
\]
since $q_m>1$. For $k\ge3$, replace $P$ by
$P/\|P\|_{\ell_2}$ and choose $\alpha\in S$ with
$|a_\alpha|=\max_{\gamma\in S}|a_\gamma|$.
Put $\Delta=\|P\|_\infty^2-1$ and $\tau=1-|a_\alpha|^2$.
As $m\ge M_k$,
\[
 0\le\Delta\le k^{1/m}-1<\theta_k,
 \qquad
 \tau\le C_k\Delta^{2^{-(k-2)}}<\kappa_k.
\]
But $k\ge3$ and $S$ is progression-free, so
Lemma~\ref{lem:noncenter} gives $\tau\ge\kappa_k$, a contradiction.
Hence $k=1$. Conversely, every nonzero monomial attains equality.
\end{proof}

\section{Proof of Theorem B(1)}\label{sec:proofB1}

\begin{proof}
Put $q_m=2m/(m+1)$ and recall $C_N$, $\kappa_N$, and $\theta_N$
from \eqref{eq:constants}--\eqref{eq:theta-M}.
Since $m\ge M_N$, we have $N^{1/m}-1<\theta_N<1$.
For monomials, equality holds. If $P$ has exactly two nonzero
coefficients, then
$\|P\|_\infty=|a_1|+|a_2|=\|P\|_{\ell_1}\ge\|P\|_{\ell_{q_m}}$.
We may therefore assume that $S:=\operatorname{supp}P$ and
$k:=|S|$ satisfy $3\le k\le N$.
Order the nonzero coefficient moduli so that
$|a_1|\ge|a_2|\ge\cdots\ge|a_k|>0$.
Suppose that
\begin{equation}\label{eq:proofB-exceptional}
 \|P\|_{\ell_{q_m}}>\|P\|_\infty.
\end{equation}
Replacing $P$ by $P/\|P\|_{\ell_2}$ leaves this inequality
and the ratio $|a_1|/|a_2|$ unchanged. Thus
$\sum_{j=1}^k|a_j|^2=1$.
Set
\[
 \Delta=\|P\|_\infty^2-1,\qquad
 \tau=1-|a_1|^2=\sum_{j=2}^k|a_j|^2.
\]
Parseval's identity, the coefficient-norm comparison, and the assumed
violation give
\begin{equation}\label{eq:proofB-delta}
 0\le\Delta<\|P\|_{\ell_{q_m}}^2-1
 \le k^{1/m}-1\le N^{1/m}-1<\theta_N<1.
\end{equation}
As $0\le\Delta<1$ and $k\le N$, Lemma~\ref{lem:concentration}
applied to $S$ gives
\begin{equation}\label{eq:proofB-tau}
\begin{aligned}
 \tau&\le C_k\Delta^{2^{-(k-2)}}\le C_N\Delta^{2^{-(N-2)}}\\
 &<C_N\theta_N^{2^{-(N-2)}}=\kappa_N.
\end{aligned}
\end{equation}
Since $|a_2|^2\le\tau$ and $|a_1|^2=1-\tau$,
\[
 \frac{|a_2|^2}{|a_1|^2}
 \le\frac{\tau}{1-\tau}<\frac{\kappa_N}{1-\kappa_N}
 =\frac1{16(N-1)^3}.
\]
Hence $|a_1|>4(N-1)^{3/2}|a_2|$, contradicting the hypothesis.
Thus \eqref{eq:theorem-B-contractivity} holds.
If the largest modulus is attained twice, then $|a_1|=|a_2|$ and
the hypothesis holds automatically. Monomials show that the
constant $1$ is optimal.
\end{proof}

\section{Proof of Theorem B(2)}\label{sec:proofB2}

The counterexamples below arise from products of Blaschke factors.
Their squared Taylor coefficient moduli define a product of geometric
distributions; the following estimate shows that a Fej\'er window of width
$m^{3/4}$ captures the fluctuations of order $\sqrt m$.

\begin{lemma}\label{lem:B2-geometric}
For every sufficiently large integer $m$, put
\[
 h=\lceil m^{3/4}\rceil,\qquad d=m-h,\qquad
 q_m=\frac{2m}{m+1},\qquad
 w_h(k)=\max\left\{1-\frac{|k|}{h+1},0\right\}\quad(k\in\mathbb Z).
\]
Then
\[
 T_m:=\sum_{\alpha\in\mathbb N_0^d}
 w_h(|\alpha|-d)^{q_m}2^{-(d+|\alpha|)q_m/2}
\]
satisfies $\lim_{m\to\infty}T_m=4$.
\end{lemma}

\begin{proof}
For integer $k$, the weight $w_h(k)$ vanishes when $|k|\ge h+1$.
Thus only indices with $d-h\le|\alpha|\le d+h$ contribute to $T_m$.
Define a measure $\nu$ on $\mathbb N_0$ by
\[
 \nu(\{j\})=2^{-(j+1)},\qquad j\ge0.
\]
Its generating function is
\[
 G(x):=\sum_{j=0}^{\infty}2^{-(j+1)}x^j=\frac1{2-x},
 \qquad |x|<2.
\]
Thus $G(1)=1$, $G'(1)=1$, and $G''(1)=2$. In particular, $\nu$ is a
probability measure with mean $1$ and variance
$G''(1)+G'(1)-G'(1)^2=2$. Under the product measure
$\nu^{\otimes d}$, the coordinates are independent, so
\begin{equation}\label{eq:B2-moments}
 \int_{\mathbb N_0^d}|\alpha|\,d\nu^{\otimes d}=d,
 \qquad
 \int_{\mathbb N_0^d}(|\alpha|-d)^2\,d\nu^{\otimes d}=2d.
\end{equation}
Since $\nu^{\otimes d}(\{\alpha\})=2^{-(d+|\alpha|)}$ and $1-q_m/2=1/(m+1)$,
\begin{equation}\label{eq:B2-integral}
\begin{aligned}
 T_m
 &=\sum_{\alpha\in\mathbb N_0^d}
       w_h(|\alpha|-d)^{q_m}2^{-(d+|\alpha|)q_m/2}\\
 &=\sum_{\alpha\in\mathbb N_0^d}
       w_h(|\alpha|-d)^{q_m}
       2^{(d+|\alpha|)/(m+1)}2^{-(d+|\alpha|)}\\
 &=\int_{\mathbb N_0^d}w_h(|\alpha|-d)^{q_m}
       2^{(d+|\alpha|)/(m+1)}\,d\nu^{\otimes d}.
\end{aligned}
\end{equation}

Let $A_m=\{\alpha\in\mathbb N_0^d:\big||\alpha|-d\big|\le m^{5/8}\}$.
On its complement, $(|\alpha|-d)^2>m^{5/4}$, and therefore
\eqref{eq:B2-moments} implies
\begin{equation}\label{eq:B2-chebyshev}
 \nu^{\otimes d}(A_m^c)
 \le\frac1{m^{5/4}}
       \int_{\mathbb N_0^d}(|\alpha|-d)^2\,d\nu^{\otimes d}
 =\frac{2d}{m^{5/4}}\le2m^{-1/4}.
\end{equation}
Since $h\ge m^{3/4}$, for all sufficiently large $m$ and every
$\alpha\in A_m$ we have $\big||\alpha|-d\big|<h+1$. Hence
\[
\begin{aligned}
 &w_h(|\alpha|-d)^{q_m}2^{(d+|\alpha|)/(m+1)}\\
 &\quad=\left(1-\frac{\big||\alpha|-d\big|}{h+1}\right)^{q_m}
    2^{(d+|\alpha|)/(m+1)}\\
 &\quad\ge\left(1-\frac{m^{5/8}}{h+1}\right)^{q_m}
    2^{(2d-m^{5/8})/(m+1)}\\
 &\quad\ge(1-m^{-1/8})^{q_m}
    2^{(2d-m^{5/8})/(m+1)}.
\end{aligned}
\]
The integrand in \eqref{eq:B2-integral} is nonnegative, so
\begin{equation}\label{eq:B2-q-lower}
\begin{aligned}
 T_m
 &=\int_{\mathbb N_0^d}w_h(|\alpha|-d)^{q_m}
    2^{(d+|\alpha|)/(m+1)}\,d\nu^{\otimes d}\\
 &\ge\int_{A_m}w_h(|\alpha|-d)^{q_m}
    2^{(d+|\alpha|)/(m+1)}\,d\nu^{\otimes d}\\
 &\ge\nu^{\otimes d}(A_m)(1-m^{-1/8})^{q_m}
    2^{(2d-m^{5/8})/(m+1)}\\
 &\ge(1-2m^{-1/4})(1-m^{-1/8})^{q_m}
    2^{(2d-m^{5/8})/(m+1)}.
\end{aligned}
\end{equation}
Since $0\le h/m\le m^{-1/4}+m^{-1}$, we have $d/m\to1$, and
\[
 \frac{2d-m^{5/8}}{m+1}
 =2-\frac{2h+m^{5/8}+2}{m+1}\longrightarrow2.
\]
Moreover, $q_m\to2$, $(1-m^{-1/8})^{q_m}\to1$, and
$1-2m^{-1/4}\to1$. Thus \eqref{eq:B2-q-lower} gives
$\liminf_{m\to\infty}T_m\ge4$.
On the other hand, $0\le w_h\le1$ and the integrand vanishes for
$|\alpha|>d+h$. Since $\nu^{\otimes d}(\mathbb N_0^d)=1$,
\[
 T_m\le 2^{(2d+h)/(m+1)}\nu^{\otimes d}(\mathbb N_0^d)
 =2^{\,2-(h+2)/(m+1)}\longrightarrow4.
\]
Therefore $T_m\to4$.
\end{proof}

\begin{proof}[Proof of Theorem B\textup{(2)}]
The construction in \cite[Proposition~6.1]{PellegrinoTeixeiraPolynomialGrowth}
combines products of Blaschke factors with Fej\'er summation. The factors
have modulus one on the circle, while the nonnegative Fej\'er kernel
selects finitely many degrees without increasing the supremum norm.
For every sufficiently large integer $m$, set
\[
 h=\lceil m^{3/4}\rceil,\qquad d=m-h,
 \qquad q_m=\frac{2m}{m+1}.
\]
Then $d>h\ge1$, $d/m\to1$, and $q_m\to2$ as $m\to\infty$.

Put $a=2^{-1/2}$ and consider the Blaschke factor
\[
\begin{aligned}
 b(w)&=\frac{a-w}{1-aw}
 =\sum_{j=0}^{\infty}c_jw^j,\qquad |w|\le1,\\
 c_0&=a,\qquad c_j=(a^2-1)a^{j-1}=-2^{-(j+1)/2}\quad(j\ge1).
\end{aligned}
\]
Moreover,
\[
 |1-aw|^2-|a-w|^2=(1-a^2)(1-|w|^2),
\]
so $|b(w)|\le1$ for $|w|\le1$, with equality when $|w|=1$.
For $z=(z_1,\ldots,z_d)\in\overline{\mathbb D}^{d}$, define
\[
 B_d(z)=\prod_{j=1}^d b(z_j)
       =\sum_{\alpha\in\mathbb N_0^d}c_\alpha z^\alpha,
 \qquad c_\alpha=\prod_{j=1}^d c_{\alpha_j}.
\]
The Taylor expansion of $B_d$ converges absolutely and uniformly on
$\overline{\mathbb D}^{d}$, since
\[
 \sum_{\alpha\in\mathbb N_0^d}|c_\alpha|
 =\left(a+(1-a^2)\sum_{j=1}^{\infty}a^{j-1}\right)^d
 =(1+2a)^d<\infty.
\]
Thus
\begin{equation}\label{eq:B2-coefficients}
 |c_\alpha|^2=2^{-(d+|\alpha|)},\qquad
 \|B_d\|_\infty=1.
\end{equation}

Define the Fej\'er kernel $F_h:\mathbb R\to[0,\infty)$ by
\[
 F_h(t)=\frac1{h+1}\left|\sum_{j=0}^{h}e^{ijt}\right|^2.
\]
For $-h\le k\le h$, there are $h+1-|k|$ pairs
$(j,\ell)\in\{0,\ldots,h\}^2$ with $j-\ell=k$. With $w_h$ as in
Lemma~\ref{lem:B2-geometric},
\begin{equation}\label{eq:B2-fejer}
 F_h(t)=\frac1{h+1}\sum_{j=0}^{h}\sum_{\ell=0}^{h}e^{i(j-\ell)t}
 =\sum_{k=-h}^{h}w_h(k)e^{ikt}.
\end{equation}
In particular, $F_h\ge0$, and integration of
\eqref{eq:B2-fejer} gives
\[
 \int_0^{2\pi}F_h(t)\,\frac{dt}{2\pi}=w_h(0)=1.
\]
Define the polynomial
\begin{equation}\label{eq:B2-Q}
 Q_m(z)=\sum_{\substack{\alpha\in\mathbb N_0^d\\
                         d-h\le|\alpha|\le d+h}}
              w_h(|\alpha|-d)c_\alpha z^\alpha.
\end{equation}
For $z\in\overline{\mathbb D}^{d}$, orthogonality and the identity
$w_h(-k)=w_h(k)$ give
\[
\begin{aligned}
 &\int_0^{2\pi}F_h(t)e^{-idt}B_d(e^{it}z)\,\frac{dt}{2\pi}\\
 &\quad=\sum_{\alpha\in\mathbb N_0^d}c_\alpha z^\alpha
           \sum_{k=-h}^h w_h(k)
           \int_0^{2\pi}e^{i(k+|\alpha|-d)t}\,\frac{dt}{2\pi}\\
 &\quad=\sum_{\alpha\in\mathbb N_0^d}
                 w_h(d-|\alpha|)c_\alpha z^\alpha
       =Q_m(z).
\end{aligned}
\]
The absolute uniform convergence of the Taylor expansion justifies termwise integration. Since $F_h\ge0$, $\int_0^{2\pi} F_h(t)\,dt/(2\pi)=1$, and
$|B_d(e^{it}z)|\le1$ for $z\in\overline{\mathbb D}^{d}$,
\[
 |Q_m(z)|
 \le\int_0^{2\pi}F_h(t)|B_d(e^{it}z)|\,\frac{dt}{2\pi}
 \le1.
\]

Define $P_m:\mathbb C^{d+1}\to\mathbb C$ by
\begin{equation}\label{eq:B2-P}
 P_m(z_0,z)=\sum_{\substack{\alpha\in\mathbb N_0^d\\
                           d-h\le|\alpha|\le d+h}}
      w_h(|\alpha|-d)c_\alpha z_0^{m-|\alpha|}z^\alpha.
\end{equation}
Since $d+h=m$, the exponent $m-|\alpha|$ is nonnegative and
every term has degree $m$. The exponent vectors
$(m-|\alpha|,\alpha)$ are distinct as $\alpha$ varies, so
homogenization does not combine coefficients. For
$(z_0,z)\in\mathbb T^{d+1}$,
\[
 P_m(z_0,z)=z_0^m Q_m(z/z_0),
\]
and conversely $P_m(1,z)=Q_m(z)$ on $\mathbb T^d$.
The maximum-modulus principle therefore gives
\begin{equation}\label{eq:B2-sup}
 \|P_m\|_\infty=\|Q_m\|_\infty\le1.
\end{equation}

For each integer $d-h\le r\le d+h$, all coefficients of
$P_m$ indexed by $\alpha\in\mathbb N_0^d$ with $|\alpha|=r$
have modulus
\[
 w_h(r-d)\,2^{-(d+r)/2}>0.
\]
Choose $r_0\in\{d-h,\ldots,d+h\}$ where this modulus is largest.
Since $r_0>0$ and $d\ge2$, the distinct monomials
\[
 z_0^{m-r_0}z_1^{r_0},\qquad z_0^{m-r_0}z_2^{r_0}
\]
both occur in $P_m$ and have coefficients of modulus
$w_h(r_0-d)2^{-(d+r_0)/2}$. Thus the largest coefficient modulus
is attained at least twice.

By the definition of the coefficients and Lemma~\ref{lem:B2-geometric},
\[
 \|P_m\|_{\ell_{q_m}}^{q_m}
 =\sum_{\alpha\in\mathbb N_0^d}
 w_h(|\alpha|-d)^{q_m}|c_\alpha|^{q_m}
 =T_m\longrightarrow4.
\]
Since $q_m\to2$,
\[
 \log\|P_m\|_{\ell_{q_m}}=\frac{\log T_m}{q_m}
 \longrightarrow\frac{\log4}{2}=\log2,
 \qquad \|P_m\|_{\ell_{q_m}}\longrightarrow2.
\]
Thus $\|P_m\|_{\ell_{q_m}}>1\ge\|P_m\|_\infty$ for all
sufficiently large $m$. Given $M\ge1$, choose such an $m\ge M$
and put $n=d+1$.
\end{proof}

\section*{Acknowledgments and Statements}

\textit{Funding.} D. Pellegrino is supported by Grants No. 406457/2023-9 (CNPq/MCTI No. 10/2023), No. 403964/2024-5 (MCTI/CNPq No. 16/2024), and No. 305807/2025-0 from the Conselho Nacional de Desenvolvimento Cient\'ifico e Tecnol\'ogico (CNPq, Brazil). E. Teixeira gratefully acknowledges support from the Grayce B. Kerr Chair at Oklahoma State University.

This work was conducted in part within the DARPA ExpMath project \emph{A Human-Centered Framework for AI-Mathematician Collaboration in Research-Level Mathematics} (Agreement No. HR0011262E029). The views and conclusions expressed here are those of the authors and should not be interpreted as representing the official policies of the Department of Defense or the U.S. Government.

\medskip

\textit{AI assistance disclosure.} The authors developed the mathematical ideas and proof strategies underlying this work and wrote the manuscript. OpenAI’s ChatGPT was used to assist with exploratory calculations, consistency checks, the organization of arguments, and editorial refinement. The authors independently reviewed and verified all AI-assisted material retained in the manuscript and assume full responsibility for the accuracy, originality, and integrity of the work.

\medskip

\textit{Competing interests.} The authors declare no competing interests.

\medskip

\textit{Data availability.} No datasets were generated or analyzed in this study.

\end{document}